\documentclass[12pt]{amsart}
\newtheorem{Theorem}{Theorem}[section]

\newtheorem{Lemma}[Theorem]{Lemma}

\theoremstyle{definition}

\theoremstyle{remark}

\begin{document}
\sloppy
\title{A new class of reconstructible graphs}
\author{Tetsuya Hosaka} 
\address{Department of Mathematics, Utsunomiya University, 
Utsunomiya, 321-8505, Japan}
\date{July 23, 2004}
\email{hosaka@cc.utsunomiya-u.ac.jp}
\subjclass[2000]{}
\thanks{Partly supported by the Grant-in-Aid for Scientific Research, 
The Ministry of Education, Culture, Sports, Science and Technology, Japan, 
(No.~15740029).}
\maketitle
\begin{abstract}
In this paper, we give a new class of reconstructible graphs.
\end{abstract}

\section{Introduction}

The purpose of this paper is to give a new class of reconstructible graphs.

For a graph $G$ and $v\in V(G)$, 
let $N(v)=\{u\in V(G)\,|\, uv\in E(G)\}$ and 
$N[v]=N(v)\cup\{v\}$.

We prove the following theorem which is an extension of a result in \cite{H}.

\begin{Theorem}
Let $G$ be a graph with $V(G)=\{v_1,\dots,v_k,\dots,v_n\}$ $(1<k<n)$.
Suppose that 
\begin{enumerate}
\item[(1)] $\bigcup_{i=1}^{k}N[v_i]=V(G)$,
\item[(2)] $N[v_1]\cap \bigcup_{i=2}^{k}N[v_i]=\emptyset$,
\item[(3)] $|d(v_i)-d(v_j)|\neq 1$ for any $i\in\{1,\dots,k\}$ 
and $j\in\{1,\dots,n\}-\{i\}$,
\item[(4)] $\{v\in V(G)\,|\, d(v)=d(v_1)\}=\{v_1\}$, and 
\item[(5)] $\{v\in V(G)\,|\, d(v)=d(v_i)\}\subset\{v_2,\dots,v_k\}$ 
for any $i\in\{2,\dots,k\}$.
\end{enumerate}
Then $G$ is reconstructible.
\end{Theorem}

\section{Proof of the theorem}

Let $G$ be a graph with $V(G)=\{v_1,\dots,v_k,\dots,v_n\}$ $(1<k<n)$.
Suppose that 
\begin{enumerate}
\item[(1)] $\bigcup_{i=1}^{k}N[v_i]=V(G)$,
\item[(2)] $N[v_1]\cap \bigcup_{i=2}^{k}N[v_i]=\emptyset$,
\item[(3)] $|d(v_i)-d(v_j)|\neq 1$ for any $i\in\{1,\dots,k\}$ 
and $j\in\{1,\dots,n\}-\{i\}$,
\item[(4)] $\{v\in V(G)\,|\, d(v)=d(v_1)\}=\{v_1\}$, and 
\item[(5)] $\{v\in V(G)\,|\, d(v)=d(v_i)\}\subset\{v_2,\dots,v_k\}$ 
for any $i\in\{2,\dots,k\}$.
\end{enumerate}
Let $G'$ be a graph with $V(G')=\{v'_1,\dots,v'_n\}$ 
such that $G-v_j \cong G'-v'_j$ for any $j\in\{1,\dots,n\}$.

Here we note that $d(v_j)=d(v'_j)$ for any $j\in\{1,\dots,n\}$.

Let $f_j:G-v_j \rightarrow G'-v'_j$ be an isomorphism 
for each $j\in\{1,\dots,n\}$.

\begin{Lemma}\label{Lemma1}
$N[v'_1]\cap \bigcup_{i=2}^{k}N[v'_i]=\emptyset$.
\end{Lemma}

\begin{proof}
Suppose that 
there exists $j\in\{1,\dots,n\}$ 
such that $v'_j\in N[v'_1]\cap \bigcup_{i=2}^{k}N[v'_i]$.
Then $v'_j\in N[v'_1]\cap N[v'_{i_0}]$ for some $i_0\in\{2,\dots,k\}$.

If $j=1$ or $j=i_0$ then $v'_{i_0}\in N[v'_1]$. 
Here we consider an isomorphism $f_{i_0}:G-v_{i_0}\rightarrow G'-v'_{i_0}$.
Since $|d(v_1)-d(v_m)|>1$ for any $m\in\{2,\dots,n\}$ by (3) and (4), 
$f_{i_0}(v_1)=v'_1$ because $d(v_1)=d(v'_1)$.
Let $d(v_1;G-v_{i_0})$ be the degree of $v_1$ in $G-v_{i_0}$.
Then 
\begin{align*}
d(v_1;G-v_{i_0})&=d(f_{i_0}(v_1);G'-v'_{i_0})=d(v'_1;G'-v'_{i_0}) \\
&=d(v'_1)-1=d(v_1)-1.
\end{align*}
Hence $v_{i_0}\in N[v_1]$. 
This contradicts (2).

Suppose that $j\not\in \{1,i_0\}$.
We consider an isomorphism $f_j:G-v_j\rightarrow G'-v'_j$.
Then $f_j(v_1)=v'_1$ by (3) and (4), 
and $f_j(v_{i_1})=v'_{i_0}$ for some $i_1\in\{2,\dots,k\}$ 
such that $d(v_{i_1})=d(v_{i_0})$ by (3) and (5).
Hence 
\begin{align*}
d(v_1;G-v_j)&=d(f_j(v_1);G'-v'_j)=d(v'_1;G'-v'_j) \\
&=d(v'_1)-1=d(v_1)-1, 
\end{align*}
and 
\begin{align*}
d(v_{i_1};G-v_j)&=d(f_j(v_{i_1});G'-v'_j)=d(v'_{i_0};G'-v'_j) \\
&=d(v'_{i_0})-1=d(v_{i_1})-1.
\end{align*}
This means that $v_j\in N[v_1]\cap N[v_{i_1}]$
which contradicts (2).

Thus $N[v'_1]\cap \bigcup_{i=2}^{k}N[v'_i]=\emptyset$.
\end{proof}

We consider an isomorphism $f_1:G-v_1\rightarrow G'-v'_1$.
By (1) and (2), 
$$f_1(V(G)- N[v_1])=f_1(\bigcup_{i=2}^{k}N[v_i]).$$
Since $f_1(\{v_2,\dots,v_k\})=\{v'_2,\dots,v'_k\}$ by (3) and (5), 
$$f_1(\bigcup_{i=2}^{k}N[v_i])=\bigcup_{i=2}^{k}N[v'_i].$$
By Lemma~\ref{Lemma1}, 
$$\bigcup_{i=2}^{k}N[v'_i] \subset V(G')- N[v'_1].$$
Thus 
$$f_1(V(G)- N[v_1])\subset V(G')- N[v'_1].$$
Here 
\begin{align*}
|f_1(V(G)- N[v_1])|&=n-d(v_1)-1 \\
&=n-d(v'_1)-1=|V(G')- N[v'_1]|.
\end{align*}
Hence $f_1(V(G)- N[v_1])= V(G')- N[v'_1]$.
This implies that $f_1(N(v_1))=N(v'_1)$.
Thus the map $f:G\rightarrow G'$ defined by 
$f(v_1)=v'_1$ and $f(v_j)=f_1(v_j)$ for each $j\in\{2,\dots,n\}$ 
is an isomorphism.
Therefore $G$ is reconstructible.

%

%

\begin{thebibliography}{1}
%
\bibitem {H}
T.~Hosaka, 
{\it A class of reconstructible graphs}, preprint.
%
\bibitem {M}
F.S.~Mulla, 
{\it A class of graphs for which the Ulam conjecture holds}, 
Discrete Math. {\bf 22} (1978), 197--198.
%
\end{thebibliography}
\end{document}